\documentclass[11pt]{amsart}

\usepackage{amsmath,amssymb,amsthm,mathtools}
\usepackage{enumitem}
\usepackage{booktabs}
\usepackage{array}
\usepackage{microtype}
\usepackage{needspace}
\usepackage[hidelinks]{hyperref}
\hypersetup{pdftitle={The Four-Model Scott Complexity Spectrum: F(4)=8},pdfauthor={Yueer Yang},pdfkeywords={Scott complexity, Ehrenfeucht theory, Vaught conjecture, Rudin--Keisler preorder, almost-prime model}}

\newtheorem{theorem}{Theorem}[section]
\newtheorem{proposition}[theorem]{Proposition}
\newtheorem{lemma}[theorem]{Lemma}
\newtheorem{corollary}[theorem]{Corollary}
\theoremstyle{definition}
\newtheorem{definition}[theorem]{Definition}
\newtheorem{remark}[theorem]{Remark}

\newcommand{\Lin}{L_{\omega_1,\omega}}
\newcommand{\Sin}[1]{\Sigma^{\mathrm{in}}_{#1}}
\newcommand{\Pin}[1]{\Pi^{\mathrm{in}}_{#1}}
\newcommand{\dSin}[1]{d\text{-}\Sigma^{\mathrm{in}}_{#1}}
\newcommand{\SC}{\operatorname{SC}}
\newcommand{\AP}{\operatorname{AP}}
\newcommand{\RK}{\operatorname{RK}}
\newcommand{\IL}{\operatorname{IL}}
\newcommand{\Th}{\operatorname{Th}}
\newcommand{\tp}{\operatorname{tp}}

\title[The Four-Model Scott Complexity Spectrum]{The Four-Model Scott Complexity Spectrum: $F(4)=8$}
\author{Yueer Yang}
\date{3 November 2026}


\makeatletter
\let\originalsetauthors\@setauthors
\renewcommand{\@setauthors}{%
  \originalsetauthors
  \vspace{-0.8em}%
  \begin{center}
    \normalfont\small Department of Mathematics, University of Maryland\\
    \footnotesize 4176 Campus Drive, William E. Kirwan Hall, College Park, MD 20742-4015, USA
  \end{center}}
\makeatother

\begin{document}

\begin{abstract}
Gonzalez, Rossegger, and Turetsky introduced a spectrum-counting function $F(n)$ for the Scott complexities realized by $\Pin{2}$ theories with exactly $n$ countable models. They proved $7\leq F(4)\leq 8$ and left the exact value open. We prove
\[
F(4)=8.
\]
More precisely, we construct a $\Pin{2}$ sentence with exactly four countable models whose optimal Scott complexities are
\[
\Pin{2},\qquad \Sin{3},\qquad \Pin{3},\qquad \dSin{3},
\]
one of each. The construction starts with a complete first-order Ehrenfeucht theory having exactly two almost-prime and two limit models, whose existence follows from a realization theorem of Sudoplatov, and then passes to a full first-order Morleyization. The proof is organized around two low-level transfer principles. First, in a small theory, a $\Sin{3}$ Scott sentence forces a model to be prime over a finite tuple; after Morleyization the converse holds. Second, if a model realizes every complete finite-tuple type over the empty set, then every $\Pin{3}$ sentence true in it is true in the countable saturated model. These two principles force the four models into the four distinct Scott-complexity classes above.
\end{abstract}

\maketitle

\section{Introduction}

For a countable structure $M$, a Scott sentence is an $\Lin$ sentence whose countable models are exactly the isomorphic copies of $M$, by Scott's theorem \cite{Scott65}. Scott complexity refines Scott rank by recording whether an optimal Scott sentence is $\Sin{\alpha}$, $\Pin{\alpha}$, or $\dSin{\alpha}$; see \cite{AGHT21,HarrisonTrainor22}. The interaction between Scott complexity and theories with finitely many countable models has recently become particularly explicit.

The finite-model setting is historically tied to Vaught's original work on denumerable models \cite{Vaught61}. Gonzalez, Rossegger, and Turetsky \cite{GRT26} define the \emph{spectrum counting function} $F\colon\omega\to\omega$ by letting $F(n)$ be the number of possible Scott-complexity spectra of $\Pin{2}$ theories with exactly $n$ countable models. They compute $F(2)=2$ and $F(3)=4$. For four models they prove
\[
7\leq F(4)\leq 8
\]
and identify the single unaccounted spectrum. In their notation, the missing possibility has one model of each of the four complexities
\begin{equation}\label{eq:missing-spectrum}
\Pin{2},\qquad \Sin{3},\qquad \Pin{3},\qquad \dSin{3}.
\end{equation}
They explicitly ask for the value of $F(4)$ in \cite[Question~3]{GRT26}. They also construct a six-model $\Pin{2}$ theory having a $\dSin{3}$ model \cite[Proposition~4.18]{GRT26}.

The purpose of this paper is to realize \eqref{eq:missing-spectrum} with exactly four models. The main structural result is the following transfer theorem.

\begin{theorem}[Four-model transfer theorem]\label{thm:transfer-intro}
Let $T_0$ be a countable complete first-order theory with exactly four countable models up to isomorphism. Suppose that exactly two of these models are prime over finite tuples. Then there is a countable definitional expansion $T_0^{\#}$ such that:
\begin{enumerate}[label=\textup{(\roman*)}]
\item $T_0^{\#}$ is a complete first-order theory with the same countable isomorphism types as $T_0$;
\item the class of models of $T_0^{\#}$ is axiomatized by a single $\Pin{2}$ sentence; and
\item the four countable models of $T_0^{\#}$ have optimal Scott complexities exactly
\[
\Pin{2},\quad \Sin{3},\quad \Pin{3},\quad \dSin{3}.
\]
\end{enumerate}
\end{theorem}

The required first-order theory $T_0$ is supplied by Sudoplatov's realization theorem for finite Rudin--Keisler preorders and distributions of limit models \cite[Theorem~4.1]{Sudoplatov06}. Taking a two-element chain and assigning two limit models to the top Rudin--Keisler class gives exactly two almost-prime models and two limit models. The standard decomposition formula, recalled in \cite[Theorem~1.10]{Sudoplatov18}, then gives four models in total.

The proof of Theorem~\ref{thm:transfer-intro} uses two facts about low-complexity Scott sentences.

First, a $\Sin{3}$ Scott sentence contains a finite existential witness followed by a $\Pin{2}$ condition. A $\Pin{2}$ condition reflects from a model to an elementary submodel containing the parameters. In a small theory, taking the prime model over the finite witness therefore shows that every model with a $\Sin{3}$ Scott sentence is almost prime. A suitable full Morleyization makes the converse true: prime models have $\Pin{2}$ Scott sentences, and models prime over finite tuples have $\Sin{3}$ Scott sentences.

Second, suppose $M$ realizes every complete finite-tuple type over $\varnothing$, and let $S$ be the countable saturated model. Any $\Pin{3}$ sentence true in $M$ transfers to $S$: one first matches the parameter tuple by its complete type, and then uses saturation of $S$ to realize the required extension type. Thus, among models realizing all $\varnothing$-types, a $\Pin{3}$ Scott sentence can only describe the saturated model.

When a small theory has exactly two almost-prime models, every nonprime model realizes all complete finite-tuple types over $\varnothing$. In the four-model situation, the models therefore arrange themselves as follows: the prime model is $\Pin{2}$; the other almost-prime model is $\Sin{3}$; the saturated model is $\Pin{3}$; and the remaining model is neither almost prime nor saturated, forcing it above both $\Sin{3}$ and $\Pin{3}$. Since only one model remains after excluding the first three, an explicit Boolean combination of their Scott sentences gives a $\dSin{3}$ Scott sentence for the fourth model.

Combining the transfer theorem with Sudoplatov's existence theorem and \cite[Proposition~4.17]{GRT26} gives the main result.

\begin{theorem}\label{thm:main}
The spectrum counting function satisfies $F(4)=8$.
\end{theorem}

\section{Preliminaries}

We work in countable languages. All structures are assumed nonempty. For the infinitary hierarchy we follow \cite{GRT26}. Thus $\Sin{0}=\Pin{0}$ is the class of finitary quantifier-free formulas. For $\alpha>0$, a formula is $\Sin{\alpha}$ if it is equivalent to a countable disjunction of formulas of the form $\exists\bar x\,\psi$, where $\psi\in\Pin{\beta}$ for some $\beta<\alpha$; $\Pin{\alpha}$ is defined dually. A $\dSin{\alpha}$ formula is a conjunction of a $\Sin{\alpha}$ formula and a $\Pin{\alpha}$ formula. We use the standard inclusions
\[
\Sin{\alpha},\Pin{\alpha}\subseteq \dSin{\alpha}
\subseteq \Sin{\alpha+1}\cap\Pin{\alpha+1}.
\]
For a countable structure $M$, $\SC(M)$ denotes the complexity of an optimal Scott sentence in the sense of \cite{AGHT21}. Infinite structures have Scott complexity at least $\Pin{2}$ \cite{AGHT21}.

If $\varphi$ is an $\Lin$ sentence with finitely many countable models and $\Gamma$ is a Scott-complexity class, write
\[
I(\varphi,\Gamma)=\bigl|\{[M]_{\cong}:M\models\varphi,\ \SC(M)=\Gamma\}\bigr|.
\]
The \emph{Scott-complexity spectrum} of $\varphi$ is the finitely supported function $\Gamma\mapsto I(\varphi,\Gamma)$. Thus $F(N)$ counts the spectra of $\Pin{2}$ sentences having exactly $N$ countable models, as in \cite[Definition~4.4]{GRT26}.

\begin{definition}
Let $T$ be a complete first-order theory. A countable model $M\models T$ is \emph{almost prime} if it is prime over a finite tuple. We write $\AP(T)$ for the set of isomorphism types of almost-prime countable models of $T$. The prime model itself corresponds to the empty tuple. Here, $M$ is prime over $\bar a$ when $(M,\bar a)$ is a prime model of its complete theory in the language with constants naming $\bar a$.
\end{definition}

The terminology agrees with the Rudin--Keisler literature used by Sudoplatov \cite{Sudoplatov06,Sudoplatov18}.

\begin{definition}
A complete first-order theory $T$ is \emph{small} if $S_n(T)$ is countable for every finite $n$.
\end{definition}

We will repeatedly use the following standard facts about small theories; see, for example, \cite[Sections~4.2--4.4]{Marker02}.

\begin{lemma}\label{lem:smallness}
Let $T$ be a countable complete first-order theory with only finitely many countable models. Then:
\begin{enumerate}[label=\textup{(\roman*)}]
\item $T$ is small;
\item for every finite tuple $\bar a$ in a model of $T$, a prime model over $\bar a$ exists; and
\item $T$ has a countable saturated model, unique up to isomorphism.
\end{enumerate}
Conclusions \textup{(ii)} and \textup{(iii)} hold for every small countable complete theory, without the assumption on the number of countable models.
\end{lemma}

\begin{proof}
Every complete finite-tuple type is realized in some countable model. A finite union of countable models realizes only countably many finite tuples, so each $S_n(T)$ is countable. If $\bar a$ is finite, then $\Th(M,\bar a)$ is again small, since an $m$-type over $\bar a$ corresponds to a complete $(m+|\bar a|)$-type over $\varnothing$ extending $\tp(\bar a)$. In a countable compact Stone space, isolated points are dense in every nonempty clopen set, so the usual atomic Henkin construction gives a countable prime model of $\Th(M,\bar a)$; its reduct is prime over $\bar a$.

For saturation, start with a countable model $M_0$. Given $M_n$, there are only countably many complete types over finite subsets of $M_n$, again by smallness over finite parameter sets. Take a countable elementary extension $M_{n+1}$ realizing all of them. The union of the resulting elementary chain is countable and realizes every type over every finite subset, hence is saturated. Uniqueness of countable saturated models follows by the standard back-and-forth argument.
\end{proof}

\section{A low-complexity full Morleyization}

We use the standard full first-order Morleyization, axiomatized by the recursive Tarskian truth clauses. The following proof records that these axioms can be chosen to have complexity $\Pin{2}$.

\begin{proposition}\label{prop:morleyization}
Let $T_0$ be a countable complete first-order theory in a countable language $L$. There is a countable definitional expansion $L^{\#}\supseteq L$ and a complete first-order $L^{\#}$-theory $T_0^{\#}$ such that:
\begin{enumerate}[label=\textup{(\roman*)}]
\item every $L$-model of $T_0$ has a unique $L^{\#}$-expansion to a model of $T_0^{\#}$;
\item $T_0^{\#}$ eliminates quantifiers;
\item elementary embeddings, primeness over finite tuples, and saturation are preserved by the expansion; and
\item there is a single $\Pin{2}$ sentence $\Theta_{T_0}$ whose models are exactly the models of $T_0^{\#}$.
\end{enumerate}
\end{proposition}

\begin{proof}
For every $L$-formula $\varphi$ and every finite list $\bar x$ of distinct variables containing its free variables, add a relation symbol $R_{\varphi,\bar x}$ of arity $|\bar x|$. We write $R_\varphi(\bar x)$ when the list is understood. This convention allows the same variable list to be used for both subformulas of a conjunction, even when some variables do not occur in one of them. We allow nullary relation symbols for sentences. If one prefers to avoid nullary symbols, they can be replaced uniformly by unary predicates forced to be constant; this does not change the argument.

Choose $\neg,\wedge,\exists$ as primitive logical operations. Add recursive defining axioms. For atomic $L$-formulas $\alpha(\bar x)$, include
\[
\forall\bar x\bigl(R_\alpha(\bar x)\leftrightarrow \alpha(\bar x)\bigr).
\]
For negation and conjunction include
\[
\forall\bar x\bigl(R_{\neg\varphi}(\bar x)\leftrightarrow \neg R_\varphi(\bar x)\bigr)
\]
and
\[
\forall\bar x\bigl(R_{\varphi\wedge\psi}(\bar x)
\leftrightarrow (R_\varphi(\bar x)\wedge R_\psi(\bar x))\bigr).
\]
For existential quantification, rename the bound variable if necessary so that $y$ is not among $\bar x$, and include
\begin{align*}
&\forall\bar x\bigl(R_{\exists y\,\psi}(\bar x)
\rightarrow \exists y\,R_\psi(\bar x,y)\bigr),\\
&\forall\bar x\forall y\bigl(R_\psi(\bar x,y)
\rightarrow R_{\exists y\,\psi}(\bar x)\bigr).
\end{align*}
Finally, for each sentence $\sigma\in T_0$, include the nullary atomic axiom $R_\sigma$.

Every axiom above is first-order $\forall\exists$ or universal. Induction on $L$-formulas shows that in every model of these axioms,
\[
R_\varphi(\bar a)\quad\Longleftrightarrow\quad \varphi(\bar a).
\]
Thus the $L$-reduct satisfies $T_0$, and the expansion is unique. Conversely, every $T_0$-model has the unique expansion obtained by interpreting $R_\varphi$ as the set defined by $\varphi$.

Given an arbitrary $L^{\#}$-formula, replace every atomic $R_\varphi$ by the corresponding $L$-formula $\varphi$. This yields an equivalent $L$-formula $\psi$. The expanded theory proves $\psi\leftrightarrow R_\psi$, so every formula is equivalent modulo $T_0^{\#}$ to an atomic formula. Hence $T_0^{\#}$ eliminates quantifiers.

Since the expansion is definitional and unique, elementary embeddings and complete first-order types are unchanged. In particular, being prime over a finite tuple and being saturated are preserved.

The axioms form a countable family of first-order $\Pi_2$ sentences. Their countable conjunction is therefore a single $\Pin{2}$ sentence, denoted $\Theta_{T_0}$, with exactly the models of $T_0^{\#}$.
\end{proof}

From now on, when $T_0$ is fixed we write $T=T_0^{\#}$ and $\Theta=\Theta_{T_0}$, and we identify an $L$-model with its unique $L^{\#}$-expansion when no confusion can arise.

\begin{remark}
We do not claim that optimal Scott complexity is invariant under this Morleyization. The expanded language is itself the language of the witness used below, and all Scott-complexity bounds are proved directly in that language. This is legitimate because $F(N)$ ranges over all countable languages.
\end{remark}

\section{\texorpdfstring{$\Sin{3}$}{Sigma-in-3} Scott sentences and almost-prime models}

The first transfer principle is a downward reflection fact.

\begin{lemma}\label{lem:pi2-down}
Let $N\preccurlyeq M$, let $\bar a\in N$, and let $\psi(\bar a)$ be a $\Pin{2}$ formula. If $M\models\psi(\bar a)$, then $N\models\psi(\bar a)$.
\end{lemma}

\begin{proof}
Put $\psi$ in the normal form
\[
\bigwedge_i\forall\bar x\,
\bigvee_j\exists\bar y\,\theta_{ij}(\bar x,\bar y,\bar a),
\]
where each $\theta_{ij}$ is finitary quantifier-free; lower-complexity pieces may be absorbed into this form. Fix $i$ and $\bar b\in N$. Since $M\models\psi(\bar a)$, for some $j$ we have
\[
M\models \exists\bar y\,\theta_{ij}(\bar b,\bar y,\bar a).
\]
This is a first-order formula with parameters from $N$. By elementarity the same existential formula holds in $N$. Since $i$ and $\bar b$ were arbitrary, $N\models\psi(\bar a)$.
\end{proof}

\begin{proposition}\label{prop:sigma3-implies-ap}
Let $T$ be a small complete first-order theory and let $M\models T$ be countable. If $M$ has a $\Sin{3}$ Scott sentence, then $M$ is almost prime.
\end{proposition}

\begin{proof}
Let
\[
\sigma_M=\bigvee_i\exists\bar x\,\psi_i(\bar x),
\qquad \psi_i\in\Pin{2},
\]
be a $\Sin{3}$ Scott sentence for $M$. Choose $i$ and a finite tuple $\bar a\in M$ with $M\models\psi_i(\bar a)$. By Lemma~\ref{lem:smallness}, a prime model $P_{\bar a}$ over $\bar a$ exists. Using its prime embedding into $M$, identify $P_{\bar a}$ with an elementary submodel of $M$ containing $\bar a$. Lemma~\ref{lem:pi2-down} gives
\[
P_{\bar a}\models\psi_i(\bar a),
\]
so $P_{\bar a}\models\sigma_M$. Since both models are countable and $\sigma_M$ is a Scott sentence, $P_{\bar a}\cong M$. Hence $M$ is prime over a finite tuple, up to isomorphism.
\end{proof}

After the Morleyization from Proposition~\ref{prop:morleyization}, the converse is also true.

\begin{proposition}\label{prop:ap-upper}
Let $T=T_0^{\#}$ be as in Proposition~\ref{prop:morleyization}.
\begin{enumerate}[label=\textup{(\roman*)}]
\item If $P$ is a countable prime model of $T$, then $P$ has a $\Pin{2}$ Scott sentence.
\item Every almost-prime countable model of $T$ has a $\Sin{3}$ Scott sentence.
\end{enumerate}
\end{proposition}

\begin{proof}
Because $T$ eliminates quantifiers, every isolated complete type is isolated by a quantifier-free formula. For each $n$, let $\mathcal I_n$ be the countable set of quantifier-free formulas that isolate complete $n$-types of $T$. Consider
\begin{equation}\label{eq:prime-scott}
\sigma_P=\Theta\wedge
\bigwedge_{n<\omega}\forall\bar x\,
\bigvee_{\theta\in\mathcal I_n}\theta(\bar x).
\end{equation}
The second conjunct says that every finite tuple realizes an isolated type, hence that the model is atomic. Formula \eqref{eq:prime-scott} is $\Pin{2}$. Every countable atomic model of a complete theory is prime, and countable prime models are unique up to isomorphism, so \eqref{eq:prime-scott} is a Scott sentence for the prime model.

Now suppose $M$ is prime over a finite tuple $\bar a$. Expand the language by constants $\bar c$ naming $\bar a$. The complete theory $\Th(M,\bar a)$ is axiomatized by $T$ together with the quantifier-free type of $\bar a$; it remains quantifier-eliminating and its countable prime model is $(M,\bar a)$. Applying the preceding construction in the language with constants yields a $\Pin{2}$ Scott sentence $\tau(\bar c)$ for $(M,\bar a)$. Replacing the constants by variables and existentially quantifying them gives
\[
\exists\bar x\,\tau(\bar x),
\]
a $\Sin{3}$ Scott sentence for $M$.
\end{proof}

Combining Propositions~\ref{prop:sigma3-implies-ap} and \ref{prop:ap-upper} gives the first key equivalence.

\begin{corollary}\label{cor:ap-sigma3}
Let $T_0$ be a small countable complete first-order theory. For a countable model $M$ of $T_0$, under the Morleyization of Proposition~\ref{prop:morleyization},
\[
M\text{ is almost prime}
\quad\Longleftrightarrow\quad
M^{\#}\text{ has a }\Sin{3}\text{ Scott sentence}.
\]
\end{corollary}

\section{Realized types and the second back-and-forth level}

We next isolate the information seen by $\Sin{2}$ and $\Pin{2}$ sentences.

\begin{lemma}\label{lem:same-types-dsigma2}
Let $M,N$ be models of the same complete first-order theory. Suppose that for every finite $n$, $M$ and $N$ realize exactly the same complete $n$-types over $\varnothing$. Then $M$ and $N$ satisfy the same $\Sin{2}$, $\Pin{2}$, and $\dSin{2}$ sentences.
\end{lemma}

\begin{proof}
It suffices to treat $\Sin{2}$. Put a $\Sin{2}$ sentence in the form
\[
\bigvee_i\exists\bar x\,\psi_i(\bar x),
\qquad \psi_i\in\Pin{1}.
\]
If $M$ satisfies it, choose $i$ and $\bar a\in M$ with $M\models\psi_i(\bar a)$. By hypothesis there is $\bar b\in N$ with
\[
\tp^M(\bar a)=\tp^N(\bar b).
\]
A $\Pin{1}$ formula is a countable conjunction of finitary universal formulas, and every such conjunct is an ordinary first-order formula in the free tuple. Hence $N\models\psi_i(\bar b)$, so $N$ satisfies the original $\Sin{2}$ sentence. The converse is symmetric. Taking negations gives agreement on $\Pin{2}$ sentences, and finite conjunctions give agreement on $\dSin{2}$ sentences.
\end{proof}

In the Morleyized language the converse is also visible: for a complete first-order type $p(\bar x)$, the sentence
\begin{equation}\label{eq:type-realization}
\exists\bar x\bigwedge_{\varphi\in p}R_\varphi(\bar x)
\end{equation}
is $\Sin{2}$ and asserts that $p$ is realized. We will only need the forward implication in Lemma~\ref{lem:same-types-dsigma2}.

\section{\texorpdfstring{$\Pin{3}$}{Pi-in-3} transfer to the saturated model}

The second transfer principle is upward rather than downward.

\begin{definition}
A model $M\models T$ is \emph{weakly saturated} if it realizes every complete finite-tuple type over $\varnothing$, i.e. every member of every $S_n(T)$.
\end{definition}

\begin{lemma}\label{lem:pi3-transfer}
Let $T$ be a complete first-order theory. Let $M\models T$ be weakly saturated, and let $S\models T$ be a countable saturated model. If $\sigma$ is a $\Pin{3}$ sentence and $M\models\sigma$, then $S\models\sigma$. Consequently, if $M$ is countable and has a $\Pin{3}$ Scott sentence, then $M\cong S$.
\end{lemma}

\begin{proof}
Write
\[
\sigma=\bigwedge_i\forall\bar x\,
\bigvee_j\exists\bar y\,\chi_{ij}(\bar x,\bar y),
\qquad \chi_{ij}\in\Pin{1}.
\]
Fix $i$ and a tuple $\bar b\in S$. Since $M$ is weakly saturated, there is $\bar a\in M$ with
\[
\tp^M(\bar a)=\tp^S(\bar b).
\]
Because $M\models\sigma$, choose $j$ and $\bar c\in M$ such that
\[
M\models\chi_{ij}(\bar a,\bar c).
\]
Transport the complete type $\tp^M(\bar c/\bar a)$ by replacing each parameter occurrence of $\bar a$ by $\bar b$. Every finite subset of the transported type is realized in $S$: its existential closure is a first-order formula in the type of $\bar a$, hence also in the type of $\bar b$. The transported type is therefore a complete type over $\bar b$ consistent with $\Th(S,\bar b)$, and saturation of $S$ gives a realization $\bar d\in S$. Consequently
\[
\tp^M(\bar a,\bar c)=\tp^S(\bar b,\bar d).
\]
Write $\chi_{ij}=\bigwedge_{k\in K_{ij}}\theta_{ijk}$, where each $\theta_{ijk}$ is a finitary universal formula. Equality of the complete types gives
\[
S\models\theta_{ijk}(\bar b,\bar d)
\qquad\text{for every }k\in K_{ij}.
\]
Hence $S\models\chi_{ij}(\bar b,\bar d)$. Thus, for every $i$ and every $\bar b\in S$, there are $j$ and $\bar d\in S$ satisfying the required disjunct, so $S\models\sigma$.

If $\sigma$ is a Scott sentence for the countable model $M$, then $S\models\sigma$ implies $S\cong M$.
\end{proof}

We also need a low upper bound for the saturated model after Morleyization.

\begin{lemma}\label{lem:sat-pi3}
Let $T=T_0^{\#}$ be the Morleyization of a small complete theory $T_0$, and let $S$ be its countable saturated model. Then $S$ has a $\Pin{3}$ Scott sentence.
\end{lemma}

\begin{proof}
Let $L$ be the language of $T_0$. For each $p\in S_n(T_0)$, set
\[
\widehat p(\bar x)=\bigwedge_{\varphi\in p}R_\varphi(\bar x).
\]
This is a $\Pin{1}$ formula. In a model of $T$, a tuple satisfies $\widehat p$ exactly when its $L$-type is $p$. For each complete type $q(\bar x,\bar y)$ of $T_0$, with $\bar x$ and $\bar y$ finite and $\bar x$ possibly empty, put $p=q\upharpoonright\bar x$ and form
\begin{equation}\label{eq:saturation-axiom}
\forall\bar x\Bigl(\neg\widehat p(\bar x)
\ \vee\ \exists\bar y\,\widehat q(\bar x,\bar y)\Bigr).
\end{equation}
The matrix is $\Sin{2}$, so \eqref{eq:saturation-axiom} is $\Pin{3}$. Since $T_0$ is small, there are only countably many such $q$. Let $\sigma_S$ be the conjunction of the single $\Pin{2}$ sentence axiomatizing $T$ from Proposition~\ref{prop:morleyization} with all sentences \eqref{eq:saturation-axiom}. Then $\sigma_S\in\Pin{3}$.

We check that a countable $M\models T$ satisfies $\sigma_S$ if and only if $M$ is saturated. Suppose first that $M$ is saturated, and fix $q$, $p$ and $\bar a\in M$ with $M\models\widehat p(\bar a)$. The set
\[
q(\bar a,\bar y)=\{\varphi(\bar a,\bar y):\varphi(\bar x,\bar y)\in q\}
\]
is a complete $L$-type over $\bar a$. It is consistent: for each finite $\Delta\subseteq q$, the formula $\exists\bar y\bigwedge\Delta(\bar x,\bar y)$ belongs to $p$, so $M$ realizes the corresponding finite subset of $q(\bar a,\bar y)$. Saturation gives $\bar b\in M$ realizing this type. Hence $M\models\widehat q(\bar a,\bar b)$, as required by \eqref{eq:saturation-axiom}.

Conversely, suppose $M\models\sigma_S$. Let $r(\bar y/\bar a)$ be any complete $L$-type over a finite tuple $\bar a\in M$, and realize it by $\bar c$ in an elementary extension $M'\succcurlyeq M$. Put
\[
q(\bar x,\bar y)=\tp^{M'\upharpoonright L}(\bar a,\bar c),
\qquad p=q\upharpoonright\bar x.
\]
Since $M\models\widehat p(\bar a)$, sentence \eqref{eq:saturation-axiom} yields $\bar b\in M$ with $M\models\widehat q(\bar a,\bar b)$. The tuples $(\bar a,\bar b)$ and $(\bar a,\bar c)$ have the same complete $L$-type, so $\bar b$ realizes $r(\bar y/\bar a)$. This proves saturation of the $L$-reduct. The definitional expansion preserves saturation by Proposition~\ref{prop:morleyization}, so $M$ is saturated as an $L^{\#}$-structure as well.

Thus the countable models of $\sigma_S$ are exactly the countable saturated models of $T$. Their uniqueness up to isomorphism shows that $\sigma_S$ is a Scott sentence for $S$.
\end{proof}

Finally, saturation and almost primeness are incompatible in the noncategorical small setting.

\begin{lemma}\label{lem:sat-not-ap}
Let $T$ be a small complete first-order theory with more than one countable model. Its countable saturated model is not almost prime.
\end{lemma}

\begin{proof}
Let $S$ be countable and saturated. Suppose toward a contradiction that $S$ is prime over a finite tuple $\bar a$. Then $(S,\bar a)$ is atomic as a model of the complete theory $\Th(S,\bar a)$. It is also saturated: every type over a finite set of parameters in the expanded language is a type over a finite set of parameters in the original language together with $\bar a$.

Hence every complete $n$-type of $\Th(S,\bar a)$ is realized in $S$ and, by atomicity, isolated. The Stone space $S_n(\Th(S,\bar a))$ is compact and discrete, therefore finite.

Now let $p\in S_n(T)$. Saturation of $S$ ensures that $p$ is realized by some $\bar b\in S$. The type $\tp(\bar b/\bar a)$ belongs to the finite set $S_n(\Th(S,\bar a))$, and forgetting the named parameters maps these finitely many types onto $S_n(T)$. Thus $S_n(T)$ is finite for every $n$. By the Ryll--Nardzewski theorem, $T$ is $\aleph_0$-categorical, contradicting the assumption that $T$ has more than one countable model.
\end{proof}

\section{Two almost-prime models force weak saturation}

The only Rudin--Keisler datum used here is the number of almost-prime isomorphism types: $|\AP(T)|=2$. The proof below uses this hypothesis and the defining embedding property of prime models. The realization theorem needed to obtain such a theory is used in Section~9.

\begin{lemma}\label{lem:two-ap-full}
Let $T$ be a small complete first-order theory with exactly two almost-prime countable models up to isomorphism. Let $P$ be the prime model and let $A$ be the other almost-prime model. Then every countable model $M\not\cong P$ is weakly saturated.
\end{lemma}

\begin{proof}
We first spell out the embedding property used below. Suppose that $\bar a$ realizes a complete type $p$ and that $P_p$ is prime over $\bar a$. If $\bar b\in M$ also realizes $p$, then $(P_p,\bar a)$ and $(M,\bar b)$ satisfy the same complete theory in the language with constants naming the tuple: each sentence in that language is an $L$-formula evaluated at the named tuple, and its truth is determined by $p$. Since $(P_p,\bar a)$ is prime for this theory, there is an elementary embedding
\[
(P_p,\bar a)\longrightarrow (M,\bar b).
\]
Its reduct is an elementary embedding $P_p\longrightarrow M$ sending $\bar a$ to $\bar b$.

Now we show that $A$ is weakly saturated. Let $p\in S_n(T)$. By Lemma~\ref{lem:smallness}, a prime model $P_p$ over a realization of $p$ exists. Since there are exactly two almost-prime isomorphism types, $P_p\cong P$ or $P_p\cong A$. In the first case $P$ realizes $p$, and because $P$ elementarily embeds into $A$, the model $A$ realizes $p$. In the second case $A$ realizes $p$ by definition. Thus $A$ realizes every finite-tuple type over $\varnothing$.

Now let $M\not\cong P$. The model $M$ is not atomic, since all countable atomic models of a complete theory are isomorphic to the prime model. Hence $M$ realizes a nonisolated complete type $q$. Choose $\bar b\in M$ realizing $q$, and let $P_q$ be prime over a realization of $q$. By the preceding argument, $P_q$ elementarily embeds into $M$, with the distinguished tuple sent to $\bar b$. Since the prime model $P$ is atomic, it cannot realize the nonisolated type $q$, so $P_q\not\cong P$. Therefore $P_q\cong A$. Thus $A$ elementarily embeds into $M$. Since $A$ is weakly saturated, so is $M$.
\end{proof}

\Needspace{12\baselineskip}
\section{The four-model transfer theorem}

We can now compute all four optimal Scott complexities.

\begin{theorem}\label{thm:transfer}
Let $T_0$ be a countable complete first-order theory with exactly four countable models up to isomorphism and exactly two almost-prime models. Let $T=T_0^{\#}$ be the Morleyization from Proposition~\ref{prop:morleyization}. Then the four countable models of $T$ have the following optimal Scott complexities, one of each:
\[
\Pin{2},\qquad \Sin{3},\qquad \Pin{3},\qquad \dSin{3}.
\]
\end{theorem}

\begin{proof}
The theory $T_0$ is small by Lemma~\ref{lem:smallness}. Let $P$ be its prime model, $A$ its other almost-prime model, and $S$ its countable saturated model. By Lemma~\ref{lem:sat-not-ap}, $S$ is not almost prime, so it is isomorphic to neither $P$ nor $A$. Let $N$ represent the fourth isomorphism type. Thus
\[
P,\qquad A,\qquad S,\qquad N
\]
are pairwise nonisomorphic. The hypothesis that there are exactly two almost-prime isomorphism types means that these are exhausted by $P$ and $A$. If $N$ were almost prime, it would therefore be isomorphic to one of them, contrary to its choice. We identify the four models with their Morleyized expansions; Proposition~\ref{prop:morleyization} preserves their isomorphism types and the property of being almost prime.

\smallskip
\noindent\emph{The prime model.}
By Proposition~\ref{prop:ap-upper}, $P$ has a $\Pin{2}$ Scott sentence. A complete theory with a finite model has only one countable isomorphism type, so all four models here are infinite. Since an infinite structure has Scott complexity at least $\Pin{2}$, we get
\[
\SC(P)=\Pin{2}.
\]

\smallskip
\noindent\emph{The nonprime almost-prime model.}
By Proposition~\ref{prop:ap-upper}, $A$ has a $\Sin{3}$ Scott sentence. By Lemma~\ref{lem:two-ap-full}, $A$ is weakly saturated, and of course the saturated model $S$ is weakly saturated. Hence $A$ and $S$ realize exactly the same complete finite-tuple types. Lemma~\ref{lem:same-types-dsigma2} shows that they satisfy the same $\dSin{2}$ sentences. Since $A\not\cong S$, no $\dSin{2}$ sentence can be a Scott sentence for $A$.

Moreover $A$ is not saturated, since the countable saturated model is unique and isomorphic to $S$. Lemma~\ref{lem:pi3-transfer} therefore shows that $A$ has no $\Pin{3}$ Scott sentence. Thus
\[
\SC(A)=\Sin{3}.
\]

\smallskip
\noindent\emph{The saturated model.}
By Lemma~\ref{lem:sat-pi3}, $S$ has a $\Pin{3}$ Scott sentence. Since $S$ is not almost prime, Proposition~\ref{prop:sigma3-implies-ap} implies that $S$ has no $\Sin{3}$ Scott sentence. Hence
\[
\SC(S)=\Pin{3}.
\]

\smallskip
\noindent\emph{The remaining model.}
The model $N$ is not almost prime, so by Proposition~\ref{prop:sigma3-implies-ap} it has no $\Sin{3}$ Scott sentence. By Lemma~\ref{lem:two-ap-full}, $N$ is weakly saturated. It is not saturated, since $N\not\cong S$. Therefore Lemma~\ref{lem:pi3-transfer} shows that $N$ has no $\Pin{3}$ Scott sentence.

It remains to give a $\dSin{3}$ upper bound. Choose Scott sentences
\[
\sigma_P\in\Pin{2},\qquad
\sigma_A\in\Sin{3},\qquad
\sigma_S\in\Pin{3}
\]
for $P,A,S$, respectively, and recall that $\Theta\in\Pin{2}$ axiomatizes $T$. Define
\begin{equation}\label{eq:N-scott}
\sigma_N=
(\neg\sigma_S\wedge\neg\sigma_P)
\wedge
(\Theta\wedge\neg\sigma_A).
\end{equation}
The first parenthesis is $\Sin{3}$ and the second is $\Pin{3}$, so $\sigma_N\in\dSin{3}$. Among countable structures satisfying $\Theta$, formula \eqref{eq:N-scott} excludes exactly the three isomorphism types $P,A,S$. Since $T$ has exactly four countable models, its remaining countable models are precisely the copies of $N$. Hence \eqref{eq:N-scott} is a Scott sentence for $N$.

Because $N$ has neither a $\Sin{3}$ nor a $\Pin{3}$ Scott sentence, no strictly lower Scott-complexity class can describe it. Therefore
\[
\SC(N)=\dSin{3}.
\]
\end{proof}

\begin{remark}\label{rem:finite-complement}
The last step uses finiteness in an essential way. Formula \eqref{eq:N-scott} is a Scott sentence because, after excluding three known isomorphism types, exactly one countable isomorphism type remains. With several remaining limit models, the same Boolean combination would define a class rather than a single isomorphism type.
\end{remark}

\section{Existence of the required first-order theory}

We now verify that the hypothesis of Theorem~\ref{thm:transfer} is realized by an Ehrenfeucht theory already supplied by the Rudin--Keisler realization theory.

Sudoplatov associates to a small complete theory $T$ the Rudin--Keisler preorder $\RK(T)$ on isomorphism types of almost-prime models. For a Rudin--Keisler equivalence class $C$, let $\IL(C)$ denote the number of isomorphism types of limit models over types represented in $C$. We only use the following two published facts.

First, Sudoplatov's realization theorem \cite[Theorem~4.1]{Sudoplatov06} says that every finite preorder $X$ with a least element $x_0$ and a greatest equivalence class $[x_1]$, together with a function
\[
f\colon X/{\sim}\to\omega
\]
satisfying the natural conditions
\[
f([x_0])=0,\qquad
f([x_1])>0\text{ if }|X|>1,\qquad
f(C)>0\text{ whenever }|C|>1,
\]
is realized by a countable small complete theory $T$ so that $\RK(T)\cong X$ and the corresponding classes have $\IL(C)=f(C)$.

Second, the characterization of finite countable-model spectra in \cite{Sudoplatov04}, restated in \cite[Theorem~1.10]{Sudoplatov18}, gives the decomposition formula
\begin{equation}\label{eq:decomp}
I(T,\aleph_0)=|\RK(T)|+\sum_{C\in \RK(T)/{\sim}}\IL(C).
\end{equation}

\begin{proposition}\label{prop:existence}
There exists a countable complete first-order theory $T_0$ with exactly four countable models, exactly two of which are almost prime.
\end{proposition}

\begin{proof}
Let $X=\{x_0<x_1\}$ be the two-element chain. Its Rudin--Keisler equivalence classes are singletons. Define
\[
f([x_0])=0,
\qquad
f([x_1])=2.
\]
The hypotheses of \cite[Theorem~4.1]{Sudoplatov06} are satisfied. Hence there is a complete theory $T_0$ with
\[
|\RK(T_0)|=2,
\qquad
\IL([x_0])=0,
\qquad
\IL([x_1])=2.
\]
Thus $T_0$ has exactly two almost-prime countable models. By \eqref{eq:decomp},
\[
I(T_0,\aleph_0)=2+0+2=4.
\]
The other two countable models are limit models and, by definition, are not prime over finite tuples.
\end{proof}

The same two-point configuration appears explicitly among the four-model Rudin--Keisler diagrams displayed in \cite[Figure~1]{Sudoplatov18}.

\section{The value of \texorpdfstring{$F(4)$}{F(4)}}

We now combine the construction with the finite-spectrum analysis of Gonzalez--Rossegger--Turetsky.

\begin{corollary}\label{cor:missing-spectrum}
There is a $\Pin{2}$ sentence with exactly four countable models whose Scott-complexity spectrum is
\[
I(\Theta,\Pin{2})=I(\Theta,\Sin{3})=I(\Theta,\Pin{3})=I(\Theta,\dSin{3})=1,
\]
and $I(\Theta,\Gamma)=0$ for every other Scott complexity $\Gamma$.
\end{corollary}

\begin{proof}
Take $T_0$ from Proposition~\ref{prop:existence}, pass to the Morleyization $T_0^{\#}$ of Proposition~\ref{prop:morleyization}, and let $\Theta$ be the single $\Pin{2}$ sentence axiomatizing it. The Morleyization preserves the four countable isomorphism types and the property of being prime over a finite tuple. Theorem~\ref{thm:transfer} computes the four optimal Scott complexities.
\end{proof}

For completeness, Proposition~4.17 of \cite{GRT26} accounts for the following seven four-model spectra below $\dSin{3}$ and proves that at most one further spectrum can occur. In the table, an entry is the number of countable models of the indicated optimal Scott complexity; each row is realized by some $\Pin{2}$ sentence.

\begin{center}
\renewcommand{\arraystretch}{1.12}
\begin{tabular}{c|ccccc}
\toprule
& $\Pin{2}$ & $\dSin{2}$ & $\Sin{3}$ & $\Pin{3}$ & $\dSin{3}$ \\
\midrule
1 & 4 & 0 & 0 & 0 & 0 \\
2 & 3 & 1 & 0 & 0 & 0 \\
3 & 2 & 2 & 0 & 0 & 0 \\
4 & 1 & 3 & 0 & 0 & 0 \\
5 & 1 & 0 & 2 & 1 & 0 \\
6 & 2 & 0 & 1 & 1 & 0 \\
7 & 1 & 1 & 1 & 1 & 0 \\
\textbf{8} & \textbf{1} & \textbf{0} & \textbf{1} & \textbf{1} & \textbf{1} \\
\bottomrule
\end{tabular}
\end{center}

The last row is exactly Corollary~\ref{cor:missing-spectrum}.

\begin{proof}[Proof of Theorem~\ref{thm:main}]
By \cite[Proposition~4.17]{GRT26},
\[
7\leq F(4)\leq 8,
\]
and the only spectrum not realized there is the last row of the table. Corollary~\ref{cor:missing-spectrum} realizes that spectrum. Hence $F(4)\geq 8$. Together with the upper bound $F(4)\leq 8$, this gives
\[
F(4)=8.\qedhere
\]
\end{proof}

\section{Further structural consequences}

Although the numerical conclusion is finite, the proof identifies a useful semantic mechanism behind the four complexities.

\begin{proposition}\label{prop:semantic-grid}
Let $T_0$ be as in Theorem~\ref{thm:transfer}, and work in its Morleyization. Then the four models are separated by the following two tests:
\begin{center}
\begin{tabular}{c|cc|c}
\toprule
Model & almost prime & saturated & Scott complexity \\
\midrule
$P$ & yes, prime & no & $\Pin{2}$ \\
$A$ & yes, nonprime & no & $\Sin{3}$ \\
$S$ & no & yes & $\Pin{3}$ \\
$N$ & no & no & $\dSin{3}$ \\
\bottomrule
\end{tabular}
\end{center}
Moreover, $A,S,N$ all realize every complete finite-tuple type over $\varnothing$.
\end{proposition}

\begin{proof}
Everything has already been proved. The almost-prime column follows from the choice of $P,A$ and Lemma~\ref{lem:sat-not-ap}. Saturation is unique, so only $S$ is saturated. Lemma~\ref{lem:two-ap-full} shows that every nonprime model is weakly saturated. The Scott complexities are those of Theorem~\ref{thm:transfer}.
\end{proof}

Thus the missing $\dSin{3}$ complexity does not arise from an arbitrary syntactic construction. It is forced by a model that simultaneously fails the two low-level semantic descriptions available to the other nonprime models: it is neither prime over a finite tuple nor saturated, while it is still indistinguishable from the saturated model at the level of which complete $\varnothing$-types are realized.

\begin{remark}
The construction also shows that a complete first-order theory, after a definitional Morleyization with a $\Pin{2}$ infinitary axiomatization, can realize the missing four-model spectrum. This is stronger than merely producing an incomplete $\Lin$ sentence with four models. It does not contradict the fact that a complete first-order theory has at most one model of Scott complexity $\Pin{2}$; see \cite[Corollary~5.3]{GRT26}. Our spectrum has exactly one such model.
\end{remark}

\begin{remark}
Gonzalez--Rossegger--Turetsky's Proposition~4.18 gives a six-model $\Pin{2}$ example containing a $\dSin{3}$ model. The present construction shows that the same top finite Scott complexity already occurs in a four-model $\Pin{2}$ theory, and in precisely the configuration required to settle $F(4)$.
\end{remark}

\section*{Acknowledgments}

I am deeply grateful to my dear family for their encouragement, and to Zhouhao, Helen, Leonid, Miss Chan, and all my friends for their kindness and support.

\medskip
\noindent\textit{Declaration of AI-assisted work.}
The main mathematical ideas and proof strategy are my own. I used OpenAI ChatGPT (in 2026, including GPT-6 during the November revision) chiefly to translate my original Chinese manuscript into English and to edit its phrasing. It also assisted with the routine technical details in Proposition~\ref{prop:morleyization} and Lemma~\ref{lem:sat-pi3}, following the mathematical approach I had chosen. In preparing this revision, I used it to help clarify proof exposition, inspect notation and LaTeX, and prepare the point-by-point response. Some of my original Chinese notes used compressed expressions influenced by Classical Chinese; these did not always translate naturally into English. I regret any obscurity in the earlier version. I remain responsible for the mathematical arguments and the final manuscript.

\end{document}